\documentclass{endm}
\usepackage{endmmacro}

\begin{document}
\begin{verbatim}\end{verbatim}\vspace{2.5cm}

\begin{frontmatter}
\title{Matrices of Forests and the Analysis of Digraphs}
\author{Pavel Chebotarev\thanksref{myemail}} and
\author{Rafig Agaev\thanksref{coemail}}
\address{Institute of Control Sciences of the Russian Academy of Sciences\\ 65 Profsoyuznaya str., Moscow 117997, Russia}
\thanks[myemail]{Email: \href{mailto:chv@lpi.ru; pchv@rambler.ru}    {\texttt{\normalshape chv@lpi.ru; pchv@rambler.ru}}}
\thanks[coemail]{Email: \href{mailto:arpo@ipu.ru; arpoye@rambler.ru} {\texttt{\normalshape arpo@ipu.ru; arpoye@rambler.ru}}}
\begin{abstract}
The matrices of spanning rooted forests are studied as a tool for analysing the structure of digraphs and
measuring their characteristics. The problems of revealing the basis bicomponents, measuring vertex
proximity, and ranking from preference relations / sports competitions are considered. It is shown that the
vertex accessibility measure based on spanning forests has a number of desirable properties. An
interpretation for the normalized matrix of out-forests in terms of information dissemination is given.
\end{abstract}
\begin{keyword}
Laplacian matrix,
spanning forest,
matrix-forest theorem,
proximity measure,
bicomponent,
ranking,
incomplete tournament,
paired comparisons

\end{keyword}
\end{frontmatter}


\newcommand\condition[2] {\medskip\noindent{\bfseries #1.\ } {\rm #2}\par\medskip}  
\newcommand\conditiont[2]{\noindent{\bfseries #1.\ }         {\rm #2}\par\medskip}  
\def\e{w}                                                
\def\G{\Gamma}                                           
\def\si{\sigma}                                          
\def\1n{1,\ldots,n}                                      
\def\0n{0,\ldots,n}                                      
\def\GG{\mathop{\mathcal G}\nolimits}                    
\def\FF{\mathop{\mathcal F}\nolimits}                    
\def\TT{\mathop{\mathcal T}\nolimits}                    
\def\PP{\mathop{\mathcal P}\nolimits}                    
\def\interca{\mathop{\scriptscriptstyle{\rm T}}\nolimits}
\def\cdc{,\ldots,}                                       
\def\ktil {\widetilde K}                                 
\def\_#1{\mathop{\hspace{-2pt}^{}_{#1}}}                 
\def\vj {\mathop{\widetilde{J}}\nolimits}                
\def\q  {\widetilde{J}}                                  
\def\J  {\widetilde{J}{}}                                
\def\x{{}}                                               
\def\z{{}}                                               
\def\beq{\begin{equation}}                               
\def\eeq{\end{equation}}                                 
\def\suml {\mathop{\sum}   \limits}                      
\def\l{\ell}                                             
\newcommand{\card}[1]{\left|#1\right|}                   
\def\sgn{\mathop{\rm sgn}\nolimits}                      
\def\adj{\mathop{\rm adj}\nolimits}                      
\def\tor{\to\hspace{-.1em}{*}\hspace{-.05em}}            
\def\rto{\hspace{-.02em}{*}\hspace{-.06em}\to}           
\def\di{d'}                                              
\def\ve{v}                                               
\def\rank{\mathop{\rm rank}\nolimits}                    
\def\a{\mathop{\alpha}\nolimits}                         
\def\id{\mathop{\rm id}\nolimits}                        
\def\od{\mathop{\rm od}\nolimits}                        
\def\tr{\mathop{\rm tr}\nolimits}                        
\def\D{{\rm\Delta}}                                      
\tolerance=800
\raggedbottom

\section{Introduction}

The matrices of routes between vertices are useful to analyse the structure of graphs and digraphs. These
matrices are the powers of the adjacency matrix. In this paper, we consider the matrices of spanning rooted
forests as an alternative tool for analyzing digraphs. We show how they can be used for finding the basis
(source) bicomponents of a digraph (Section~\ref{stru}), for measuring vertex proximity
(Section~\ref{dosti}), and for ranking on the base of preference relations / sports competitions
(Section~\ref{leader}). In the initial sections, we introduce the necessary notation (Section~\ref{Notatio}),
present recurrent formulae for calculating the ``forest matrices'' (Section~\ref{sec_constr}), and list some
properties of spanning rooted forests and forest matrices (Section~\ref{sect_prop}).

Three features that distinguish the matrices of forests from the matrices of routes are notable. First, all
column sums (or row sums) of the forest matrices are the same, therefore, these matrices can be considered as
matrices of {\em relative\/} accessibility. Second, there are matrices of ``out-forests'' and matrices of
``in-forests'', enabling one to distinguish ``out-accessibility'' from ``in-accessibility'', which is
intuitively defensible. Third, the total weights of maximum spanning forests are closely related to the
Ces\'aro limiting probabilities of Markov chains determined by the digraph under consideration.

\section{Notation and simple facts}
\label{Notatio}

\subsection{Weighted digraphs, components, and bases}

Suppose that $\G$ is a weighted digraph without loops, $V(\G)=\{\1n\},$
$n>1,$ is its set of vertices and $E(\G)$ the set of arcs.
Let $W=(\e\_{ij})$ be the matrix of arc weights. Its entry $\e\_{ij}$ is zero if there is no arc from vertex
$i$ to vertex~$j$ in~$\G$; otherwise $\e\_{ij}$ is strictly positive. In what follows, $\G$ is fixed, unless
otherwise specified. If $\G'$ is a subgraph of $\G$, then the weight of $\G'$, $\e(\G')$, is the product of
the weights of all its arcs; if $E(\G')=\varnothing$, then $\e(\G')=1$ by definition. The weight of a
nonempty set of digraphs $\GG$ is
\beq
\label{set_weight}
\e(\GG)=\suml_{H\in\GG}\e(H);\quad \e(\varnothing)=0.
\eeq

{\it A spanning subgraph\/} of $\G$ is a subgraph with vertex set $V(\G)$. The {\it indegree\/} $\id(\ve)$
and {\it outdegree\/} $\od(\ve)$ of a vertex $\ve$ are the number of arcs that come {\em in} $\ve$ and {\em
out of} $\ve$, respectively. A vertex $\ve$ is called a {\it source\/} if $\id(\ve)=0$. A vertex $\ve$ is
{\it isolated\/} if $\id(\ve)=\od(\ve)=0$. A {\it route\/} ({\it semiroute}) is an alternating sequence of
vertices and arcs $\ve\_0,e\_1,$ $\ve\_1\cdc e\_k, \ve\_k$ with every arc $e\_i$ being $(\ve\_{i-1},\ve\_i)$
(resp., either $(\ve\_{i-1},\ve\_i)$ or $(\ve\_i,\ve\_{i-1})$). A {\it path\/} is a route with distinct
vertices. A~{\it circuit\/} is a route with $\ve\_0=\ve\_k$, the other vertices being distinct and different
from $\ve\_0$. Vertex $\ve$ {\it is reachable\/} from vertex $z$ in $\G$ if $\ve=z$ or $\G$ contains a path
from $z$ to~$\ve$.

A digraph is {\it strongly connected\/} (or {\it strong}) if all of its vertices are mutually reachable and
{\it weakly connected\/} if any two different vertices are connected by a semiroute. Any maximal strongly
connected (weakly connected) subgraph of $\G$ is a {\it strong component}, or a {\it bicomponent} (resp., a
{\it weak component}) of~$\G$. Let $\G_1\cdc\G_r$ be all the strong components of~$\G$. The {\it
condensation\/} (or {\it factorgraph}, or {\it leaf composition}, or {\it Hertz graph}) $\G^{\circ}$ of
digraph $\G$ is the digraph with vertex set $\{\G_1\cdc\G_r\}$, where arc $(\G_i,\G_j)$ belongs to
$E(\G^{\circ})$ iff $E(\G)$ contains at least one arc from a vertex of $\G_i$ to a vertex of~$\G_j$. The
condensation of any digraph $\G$ obviously contains no circuits.

A {\em vertex basis\/} of a digraph $\G$ is any minimal (by inclusion) collection of vertices such that every
vertex of $\G$ is reachable from at least one vertex of the collection. If a digraph does not contain
circuits, then its vertex basis is obviously unique and coincides with the set of all
sources~\cite{Harary69,Zykov69}. That is why the bicomponents of $\G$ that correspond to the sources of
$\G^{\circ}$ are called the {\it basis bicomponents\/}~\cite{Zykov69} or {\it source bicomponents\/} of~$\G$.
In this paper, the term {\it source knot of\/}~$\G$ will stand for the set of vertices of any source
bicomponent of~$\G.$ In~\cite{FiedlerSedlacek58}, source knots are called {\em W-bases}.

The following statement \cite{Harary69,Zykov69} characterizes all the vertex bases of a digraph.

\begin{proposition}
\label{proZy} A set $U\subseteq V(\G)$ is a vertex basis of\/ $\G$ if and only if\/ $U$ contains exactly one
vertex from each source knot of\/ $\G$ and no other vertices.
\end{proposition}

\label{sec2}

Schwartz \cite{Schwartz86} referred to the source knots of a digraph as {\em minimum $P$-undominated sets}.
According to his {\em Generalized Optimal Choice Axiom\/} (GOCHA), if a digraph represents a preference
relation on a set of alternatives, then the {\it choice\/} should be the union of its minimum $P$-undominated
sets.\footnote{This union is also called the {\em top cycle\/} and the {\em strong basis\/} of~$\G$.} This
choice is interpreted as the set of ``best'' alternatives. A review of choice rules of this kind can be found
in \cite{Vol'skii88}; for ``fuzzy'' extensions, see~\cite{Roubens96FSS}.

\subsection{Matrices of forests}

A {\it diverging tree\/} is a weakly connected digraph in which one vertex (called the {\it root}) has
indegree zero and the remaining vertices have indegree one. A~diverging tree is said to {\em diverge from\/}
its root. Spanning diverging trees are sometimes called {\it out-arborescences}. A~{\it diverging forest\/}
(or {\it diverging branching}) is a digraph all of whose weak components are diverging trees. The roots of
these trees are called the roots of the diverging forest. A~{\em converging tree\/} ({\em converging forest})
is a digraph that can be obtained from a diverging tree (resp., diverging forest) by the reversal of all
arcs. The roots of a converging forest are its vertices that have outdegree zero. In what follows, spanning
diverging forests in $\G$ will be called {\it out-forests\/} of $\G$; spanning converging forests in $\G$
will be called {\it in-forests\/} of $\G$.

\begin{definition}
\label{De-Max}
An out-forest $F$ of a digraph $\G$ is called a {\em maximum out-forest\/} of $\G$ if $\G$ has
no out-forest with a greater number of arcs than in~$F$.
\end{definition}

It is easily seen that every maximum out-forest of $\G$ has the minimum possible number of diverging trees;
this number will be called the {\it out-forest dimension\/} of $\G$ and denoted by~$\di$. It can be easily
shown that the number of arcs in any maximum out-forest is $n-\di$; in general, the number of weak components
in a forest with $k$ arcs is $n-k$.

By $\FF^{\rto}(\G)=\FF^{\rto}$ and $\FF^{\rto}_k(\G)=\FF^{\rto}_k$ we denote the set of all out-forests of
$\G$ and the set of all out-forests of $\G$ with $k$ arcs, respectively; $\FF^{i\rto j}_k$ will designate the
set of all out-forests with $k$ arcs where $j$ belongs to a tree diverging from~$i$; $\FF^{i\rto
j}=\bigcup_{k=0}^{n-\di}\FF^{i\rto j}_k$ is the set of such out-forests with all possible numbers of arcs.
The notation like $\FF^{\rto}_{(k)}$ will be used for sets of out-forests that consist of $k$ trees, so
$\FF^{\rto}_{(k)}=\FF^{\rto}_{n-k},\;k=\1n.$ Thus, the ${*}\hspace{-.4em}\to$ sign relates to out-forests;
the corresponding notation with $\to\hspace{-.3em}{*}\,$, such as $\FF^{\tor},$ relates to in-forests, i.e.,
$\ast$ images the root(s).

Let

\beq
\label{sik}
\si\_k=\e(\FF^{\rto}_k),\quad k=0,1,\ldots;\quad
\si=\e(\FF^{\rto})=\suml_{k=0}^{n-\di}\si\_k.
\label{si}
\eeq

By (\ref{sik}) and (\ref{set_weight}), $\si\_k=0$ whenever $k>n-\di;$ $\si\_0=1.$

We also introduce the parametric value
\beq
\label{sitau}
\si(\tau)
=\suml_{k=0}^{n-\di}\e(\FF^{\rto}_k)\,\tau^k
=\suml_{k=0}^{n-\di}\si\_k\tau^k,\quad \tau>0,
\eeq
which is the total weight of out-forests in $\G$ provided that all arc weights are multiplied by~$\tau$.

Consider the {\em matrices $Q\_k=(q_{ij}^k),\; k=0,1,\ldots,$ of out-forests of\/ $\G$ with $k$ arcs}: the
entries of $Q\_k$ are
\beq
\label{qijk}
q_{ij}^k=\e(\FF_k^{i\rto j}).
\eeq

By (\ref{qijk}) and (\ref{set_weight}), $Q\_k=0$ whenever $k>n-\di;$
$Q\_0=I.$

The {\em matrix of all out-forests\/} is
\beq
\label{qij}
Q=(q\_{ij})=\suml_{k=0}^{n-\di}Q\_k
\mbox{~~with entries~~}
q\_{ij}=\e(\FF^{i\rto j}).
\eeq

We will also consider the {\em normalized matrices of out-forests}:
\beq
\label{Jk}
\label{J}
J\_k=\si_k^{-1}Q\_k,\quad k=\0n-\di;
\quad
J=(J\_{ij})=\si^{-1}Q
\eeq
and the parametric matrices
\beq
\label{Qtau}
\label{Jtau}
Q(\tau)=\suml_{k=0}^{n-\di} Q\_k\tau^k
\mbox{  and  }
J(\tau)=\si^{-1}(\tau)\,Q(\tau),
\quad \tau>0,
\eeq
where $\si\_k,$ $\si,$ and $\si(\tau)$ are defined by
(\ref{sik}) and (\ref{sitau}).

The {\em normalized matrix of maximum out-forests\/} $J\_{n-\di}\/$
will also be denoted by~$\J=(\J\_{ij})$:
\beq
\label{Jbar}
\J=J\_{n-\di}.
\eeq

\section{Counting forests by means of linear algebra}
\label{sec_constr}
\label{Sec_Poli}

An algorithmic description of maximum out-forests was given in~\cite{AgaChe00}. Another algorithm for the
enumeration of out-forests can be obtained by adding a fictitious vertex $0$ along with the arcs $(0,v)$ for
all vertices $v$ of $\G$ and enumerating the diverging trees in the supplemented digraph.

The {\it column Laplacian\/} matrix of $\G$ is the $n\times n$ matrix $L=L(\G)=(\l\_{ij})$ with entries
$\l\_{ij}=-\e\_{ij}$ whenever $j\ne i$ and $\l\_{ii}=-\suml_{k\ne i}\l\_{ki}$, $i,j=\1n$. This matrix has
zero column sums; in \cite{CheAga02a+} we denoted it $L'$; now, for simplicity, designation $L$ is used. A
{\it row Laplacian\/} matrix differs from the column Laplacian matrix by the diagonal only: its diagonal is
such that the row sums are zero. The row and column Laplacian matrices are singular M-matrices (see, e.g.,
\cite[p.~258]{CheAga02a+}). Their index is 1~\cite[Proposition~12]{CheAga02a+}. The spectra of the row
Laplacian matrices were studied in~\cite{AgaChe05LAA}.

The following theorem provides a method to calculate the forest matrices by means of linear algebra.

\begin{theorem}
\label{pro.allk}
$\!\!\!${\rm~\cite{CheAga02a+}.}
${\displaystyle\; Q_{k+1}=\si\_{k+1}\!I-LQ_{k};\;\;\;
\si_{k+1}=\frac{\tr(LQ\_k)}{k+1},
\;\;\;k=0,1,\ldots.}$
\end{theorem}

Hence,
\beq
Q_{k+1}=\frac{\tr(LQ\_k)}{k+1}I-LQ_{k},\;\;\;k=0,1,\ldots.
\eeq

Consider the matrices $L\_{k}\stackrel{{\rm def}}{=}\si\_{k}\!I-Q\_{k},$ $k=0,1,\ldots.$ Obviously, $L\_{k}$
is the Laplacian matrix $L(\G^{k})$ of the {\em digraphs $\G^{k}$ of out-forests of\/~$\G$}: the arc weights
in $\G^{k}$ are the off-diagonal entries of $Q\_{k}.$ Then, by Theorem~\ref{pro.allk}, we have

\begin{corollary}
$L\_{k+1}=LQ\_k$ and $\,\tr(L\_k)=k\si\_k,\;k=0,1,\ldots.$
\end{corollary}
This implies the following recurrent formula for $L\_{k+1}$:

\begin{corollary}
${\displaystyle L\_{k+1}=L\left(\frac{\tr(L\_k)}{k}I-L\_k\right),\;k=1,2,\ldots.}$
\end{corollary}

\section{Properties of the forest matrices}
\label{sect_prop}

A number of results on the forest matrices are presented in~\cite{CheAga02a+}. Some of them are collected in
the following theorem.

\begin{theorem}
\label{th4}
\label{sumsum7}
\label{teo.allk}
$\!\!\!${\rm~\cite{CheAga02a+}.}
{\rm 1.} Matrices $J\_k,\;k=\0n-\di,$ $J,$ and $J(\tau)$ are column
         stochastic.\\
{\rm 2.} For any $\tau>0,$ $Q(\tau)=\adj(I+\tau L)$ and
         $\si(\tau)=\det(I+\tau L)$ hold$,$ whence$,$
         $J(\tau)=(I+\tau L)^{-1}.$\\
{\rm 3.} $L\q=\q L=0$.\\
{\rm 4.} $\vj$ is idempotent$:$ $\;\q^{2}=\q.$\\
{\rm 5.} $\J
         =\lim_{\tau\to\infty}J(\tau)
         =\lim_{\tau\to\infty} (I+\tau\,L)^{-1}.$\\
{\rm 6.} $\rank(\J)=\di;\;\rank(L)=n-\di$.\\
{\rm 7.} ${Q\_{k}
         =\sum_{i=0}^k\si\_{k-i}(-L)^i,\;\;\; k=0,1,\ldots.}$\\
{\rm 8.} $\J$ is the eigenprojection of $L$.
\end{theorem}

Item 2 of Theorem~\ref{th4} is a parametric version of the {\em matrix-forest theorem\/}~\cite{CheSha97}.

To formulate the topological properties of the matrix $\J$, the following notation is needed.

Let $\ktil=\bigcup_i K\_i$, where $K\_i$ are all the source knots of~$\G$; let $K_i^{+}$ be the set of all
vertices reachable from $K\_i$ and unreachable from the other source knots. For any $k\in\ktil$, $K(k)$ will
designate the source knot that contains~$k$. For any source knot $K$ of $\G,$ denote by $\G_K$ the
restriction of $\G$ to $K$ and by $\G_{-K}$ the subgraph with vertex set $V(\G)$ and arc set $E(\G)\setminus
E(\G_K)$. For a fixed $K$, $\TT$ will designate the set of all spanning diverging trees of $\G_K$, and $\PP$
the set of all maximum out-forests of $\G_{-K}$. By $\TT^k,$ $k\in K,$ we denote the subset of $\TT$
consisting of all trees that diverge from $k$, and by $\PP^{K\rto j},$ $j\in V(\G),$ the set of all maximum
out-forests of $\G_{-K}$ such that $j$ is reachable from some vertex that belongs to~$K$ in these forests.
$\q_{k\bullet}$ is the $k$th row of~$\q$.

\begin{theorem} 
\label{th2}
$\!\!\!${\rm~\cite{AgaChe00}.}
{Let $K$ be a source knot in~$\G$. Then the following statements
hold.\\
{\rm 1.}~$\q_{ij}\ne 0\;\Leftrightarrow\; (i\in\ktil$
         and $j$ is reachable from $i$ in~$\G).$\\
{\rm 2.}~Let $k\!\in\!K.\!$ For any $j\in V(\G),$
         $\q_{kj}=\e(\TT^k)\e(\PP^{K\rto j})\slash
         \e(\FF^{\rto}_{(\di)}).\!$ Furthermore$,$ if $j\in K^+,$ then
         $\q_{kj}=\q_{kk}=\e(\TT^k)\slash \e(\TT).$\\
{\rm 3.}~$\suml_{k\in K}\q_{kk}=1.$ In particular$,$
         if $k$ is a source$,$ then $\q_{kk}=1.$\\
{\rm 4.}~For any $k\_1,k\_2\in K$, $\q_{k\_2\!\bullet}=(\e(\TT^{k\_2})
         \slash\e(\TT^{k\_1}))\q_{k\_1\!\bullet}$ holds$,$ i.e.$,$ the
         rows $k\_1$ and $k\_2$ of $\q$ are proportional.}
\end{theorem}

We say that a weighted digraph $\G$ and a finite homogeneous Markov chain with transition probability matrix
$P$ {\em inversely correspond\/} to each other if \beq \label{G_M_cor} I-P=\a L^{\interca}, \eeq where $\a$
is any nonzero real number.

If a Markov chain inversely corresponds to $\G,$ then the probability of transition from $j$ to $i\ne j$ is
proportional to the weight of arc $(i,j)$ in $\G$ and is $0$ if $E(\G)$ does not contain $(i,j).$ We consider
such an {\em inverse\/} correspondence in order to model preference digraphs in Section~\ref{leader}: in this
case, the transitions in the Markov chain are performed from ``worse'' objects to ``better'' ones, so the
Markov chain stochastically ``searches the leaders.''

\begin{theorem}
\label{M} For any finite Markov chain$,$ its matrix of Ces\'aro limiting probabilities coincides with the
matrix $\J$ of any digraph inversely corresponding to this Markov chain.
\end{theorem}

Theorem~\ref{M} follows from the {\em Markov chain tree theorem\/} \cite{LeightonRivest83,LeightonRivest86},
which, in turn, can be immediately proved using item~8 of Theorem~\ref{th4} and a result
of~\cite{Rothblum76a} (see~\cite{CheAga02a+}). Another proof of Theorem~\ref{M} can be found in
\cite{chebotarev02spanning}. A review on forest representations of Markov chain probabilities is given
in~\cite{Che04MCTT}. For an interpretation of $J(\tau)$ in terms of Markov chains we refer
to~\cite{AgaChe01}.

\section{Detecting the source knots of a digraph}
\label{stru}

In this section, we show that entry $ij$ of $(I+\tau L)^{-1}$, where $\tau>0$, is nonzero if and only if $j$
is reachable from $i$ in $\G$ and that $\J$ points out the source knots of~$\G$ and the vertices reachable
from each of them.

The {\em reachability matrix\/} of a digraph is the matrix
$R=(r\_{ij})_{n\times n}$ with entries
\[
r\_{ij}=
\cases{
1, &if $j$ is reachable from $i$,\cr
0, &otherwise.
}
\]

It follows from the definition (\ref{Jtau}) of $J(\tau)$ that for every $\tau>0$
\beq
R=\sgn(J(\tau)),
\eeq
where the signum function operates entrywise.

Recall that, by item~2 of Theorem~\ref{teo.allk}, $J(\tau)$ can be calculated as follows:
$J(\tau)=(I+\tau L)^{-1}.$ So we obtain

\begin{proposition}
\label{reach} For every $\tau>0,$ $R=\sgn\left((I+\tau L)^{-1}\right)$.
\end{proposition}

An algebraic way to recover the bicomponents of a digraph is to calculate the {\it mutual reachability
matrix}, which is the Hadamard (entrywise) product of the reachability matrix and its transpose. Note that
the standard algebraic means of finding the reachability matrix of a digraph is to compute $(I+A)^{n-1}$,
where $A$ is the adjacency matrix,  or to successively calculate the power matrices $(I+A)^k$ until the
positions of nonzero entries stabilize; then, in both cases, the nonzero entries of the resulting matrix
should be replaced by ones~\cite{Zykov69}. The matrices $(I-\a A)^{-1}$ with sufficiently small $\a>0$ can
also be used for this end.

The following result does not provide an effective way to find the reachability matrix, but it contributes to
the understanding of the nature of the forest matrices.

By (\ref{Jtau}), $J(\tau)=\si^{-1}(\tau)\suml_{k=\di}^{n} Q\_{(k)}\tau^{n-k}$, where $Q\_{(k)}\stackrel{{\rm
def}}{=}Q\_{n-k}.$ It turns out that all information about the digraph reachability is accumulated in
$Q\_{(\di)}$ and $Q\_{(\di+1)}$. This follows from

\begin{proposition}
\label{p1.2} For any $i,j\in V(\G)$ and any path from $i$ to $j$ in\/ $\G,$ there exists an out-forest in
$\FF^{i\rto j}_{(\di)}\cup\FF^{i\rto j}_{(\di+1)}$ that contains this path.
\end{proposition}

\begin{proof}
{For the given path and any maximum out-forest in $\G$, consider their join and remove all arcs of the
out-forest that come in the vertices of the path but do not belong to the path. The resulting subgraph
contains neither circuits nor vertices $\ve$ with $\id(\ve)>1,$ i.e., it is an out-forest. It is rooted at
$i$ and contains at least $n-\di-1$ arcs, including the arcs of the given path. Hence, it belongs to
$\FF^{i\rto j}_{(\di)} \cup\FF^{i\rto j}_{(\di+1)}.$}
\end{proof}

This implies
\begin{corollary}
\label{p1.3}
$R=\sgn\!\left(J\_{(\di)}+J\_{(\di+1)}\right),$ where
$J\_{(k)}\stackrel{{\rm def}}{=}J\_{n-k}.$
\end{corollary}

In some cases, the main goal is to find the source knots and the vertices reachable from each of them. For
example, as was noted in Section~\ref{sec2}, this is the case in choice theory if the Generalized Optimal
Choice Axiom (GOCHA) is adopted. Then the union $\ktil=\bigcup^{\di}_{i=1}K_i$ of the source knots is the set
of ``best'' alternatives chosen on the base of a preference relation or a digraph of
preferences~\cite{Schwartz86} (cf.~\cite{Vol'skii88,Laslier97,Roubens96FSS}).

It turns out that $\q$ immediately reveals $\ktil$.

\begin{definition}
\label{De2}
The {\em top reachability matrix\/} of $\G$ is the matrix
$\widehat R=(\widehat r\_{ij})_{n\times n}$ with entries
\beq
\label{reachM}
\widehat r\_{ij}=\cases{
1, &if $i$ belongs to a source knot of $\G$ and $j$ is reachable
    from $i$,\cr
0, &otherwise.
}
\eeq
\end{definition}

It follows from item~1 of Theorem~\ref{th2} that

\begin{proposition}
\label{TR}
$\widehat R=\sgn(\J).$
\end{proposition}

The source knots of $\G$ can be disclosed by computing the {\em mutual top reachability matrix}, which is the
Hadamard product of $\widehat R$ and ${\widehat R}^{\interca}$. It follows from Proposition~\ref{TR} that
\begin{proposition}
\label{p_for_re}
$\!$Vertices $i$ and $j$ belong to the same source knot iff ${\J\_{i\!j}\J\_{ji}\!\ne\!0.}$
\end{proposition}

If $\widehat R$ is found by of approximate calculation based on item~5 of Theorem~\ref{th4}, the following
statement can be of help. Recall that, by item~2 of Theorem~\ref{th4}, $\si=\det(I+L).$

\begin{proposition}
\label{approx}
Let\/ $\G$ be a digraph with all arc weights~1. Then
\beq
\widehat r\_{ij}
=\cases{
0, &if  $J\_{ij}(\si^2)<\si^{-1},$\cr
1, &otherwise,
}
\quad i,j=\1n,
\eeq
where $J\_{ij}(\tau)$ is the $ij$-entry of $J(\tau)=(I+\tau L)^{-1}.$
\end{proposition}

Proposition~\ref{approx} is formulated for unweighted digraphs, since the reachability relation does not
depend on the weights of arcs.

\begin{proof}
{If $\di=n,$ the claim is obvious. Let $\di<n.$ Then, by (\ref{sitau}) and
(\ref{Jtau}), $J(\tau)=B(\tau)+C(\tau)$ holds, where
$
B(\tau)=(b_{ij})={\si^{-1}(\tau)}\sum_{k=0}^{n-\di-1}\,\tau^k\,Q_k
$
and
$
C(\tau)=(c_{ij})=\si^{-1}(\tau)\,\tau^{n-\di}\,Q_{n-\di}.
$
By item~5 of Theorem~\ref{th4}, we have
$B(\tau)\to 0,\;$ $J(\tau)\to\J$ and $C(\tau)\to\J$ as $\tau\to\infty.$

Let $\tau=\si^2.$ Since the weights of all arcs are 1,
$\si\_k\ge1,$ $k=\0n-\di,$ and $\tau\ge1$ hold.

Let $\widehat r\_{ij}=1.$ Then, by Proposition~\ref{TR}, $\J\_{ij}>0,$
hence, $q_{ij}^{n-\di}\ge1$, where $\left(q_{ij}^k\right)=Q_k,\;k=\0n-\di.$ Therefore,
$
J\_{ij}(\tau)
>c\_{ij}(\tau)
\ge\si^{-1}(\tau)\tau^{n-\di}
=(\sum_{k=0}^{n-\di}\si\_k\tau^k)^{-1}\tau^{n-\di}
>(\si\tau^{n-\di})^{-1}\tau^{n-\di}=\si^{-1}.
$

Let now $\widehat r\_{ij}=0.$ Then $\J\_{ij}=0,$ hence, $c\_{ij}(\tau)=0$
and
$
J\_{ij}(\tau)
=b\_{ij}(\tau)
\le(\sum_{k=0}^{n-\di}\si\_k\tau^k)^{-1}
    \sum_{k=0}^{n-\di-1}\tau^{n-\di-1}q_{ij}^k
<(\tau^{n-\di})^{-1}\tau^{n-\di-1}\si
=\si/\tau
=\si^{-1}.
$
}
\end{proof}

\section{Forest based accessibility measures}
\label{dosti}

Formally, by an {\em accessibility measure\/} for digraph vertices we mean any function that assigns a matrix
$P=(p\_{ij})\_{n\times n}$ to every weighted digraph $\G,$ where $n=\card{V(\G)}.$ Entry $p\_{ij}$ is
interpreted as the accessibility (or connectivity, relatedness, proximity, etc.) of $j$ from $i.$

Consider the accessibility measures $P^{{\rm out}}_{\tau}=J(\tau),$ where $J(\tau)$ is defined by
(\ref{Jtau}), and
$P^{{\rm in}}_{\tau}=(p^{{\rm in}}_{ij})$ with
$p^{{\rm in}}_{ij}=\e(\FF^{i\tor j}(\tau))/\e(\FF^{\tor}(\tau)),$ where
$\FF^{i\tor j}(\tau)$ and
$\FF^{\tor}(\tau)$ are, respectively, the
$\FF^{i\tor j}$ and
$\FF^{\tor}$ for the digraph $\G(\tau)$ obtained from $\G$ by the multiplication of all arc weights
by~$\tau.$ Parameter $\tau$ specifies the relative weight of short and long ties in~$\G$.

\begin{definition}
\label{dua}
Accessibility measures $P^{(1)}$ and $P^{(2)}$ are {\em dual\/} if for
every $\G$ and every $i,j\in V(\G),\;$
$p^{(1)}_{ij}(\G)=p^{(2)}_{ji}(\G'),$ where $\G'$ is obtained from $\G$
by the reversal of all arcs (preserving their weights).
\end{definition}

The following proposition results from the fact that the reversal of all arcs in $\G$ transforms all
out-forests into in-forests and vice versa.

\begin{proposition}
\label{pr_dual}
For every $\tau>0,$ the measures $P^{{\rm out}}_{\tau}$ and
                                 $P^{{\rm  in}}_{\tau}$ are dual.
\end{proposition}

What is the difference in interpretation between
$P^{{\rm out}}_{\tau}$ and
$P^{{\rm in}}_{\tau}$? A partial answer is as follows.
$P^{{\rm out}}_{\tau}$ can be interpreted as the relative weight of
$i\to j$ connections among the out-connections of $i,$ whereas
$P^{{\rm in}}_{\tau}$ is the relative weight of
$i\to j$ connections among the in-connections of $j.$ Naturally, these
relative weights need not coincide. For example, a connection between
an average man and a celebrity is usually more important for the
average man. This example demonstrates that self-duality is not an
imperative requirement to accessibility measures.
The properties of several self-dual measures have been studied in~\cite{CheSha98}.

The following conditions proposed in part in \cite{CheSha98} can be considered as desirable properties of vertex
accessibility measures.

\condition{Nonnegativity}
{$p\_{ij}\ge0,\;\:i,j\in V(\G)$.}

\conditiont{Reachability condition}
{For any $i,j\in V(\G),\;$ ($p\_{ij}=0 \Leftrightarrow j$ is unreachable from~$i$).}

\conditiont{Self-accessibility condition}
{$\!\!\!$For any distinct $i,j\in V(\G),$ (A)~$p\_{ii}>p\_{ij}$ and (B)~$p\_{ii}>p\_{ji}$ hold.}

\conditiont{Triangle inequalities for proximities} {For any $i,k,t\in V(\G)$,
(A)~$p\_{ki}-p\_{ti}$ $\le p\_{kk}-p\_{tk}$ and
(B)~$p\_{ik}-p\_{it}$ $\le p\_{kk}-p\_{kt}$ hold.
}

The triangle inequalities for proximities a counterparts of the ordinary triangle inequality which
characterizes distances (cf.~\cite{CheSha98a}).

Let $k,i,t\in V(\G)$. We say that $k$ {\em mediates between $i$\/ and $t$\/} if\/ $\G$ contains a path from $i$
to $t,$ $i\ne k\ne t,$ and every path from $i$ to $t$ includes~$k.$

\condition{Transit property}
{If $k$ mediates between $i$\/ and $t,$ then (A)~$p\_{ik}>p\_{it}$ and (B)~$p\_{kt}>p\_{it}.$ }

\conditiont{Monotonicity} {Suppose that the weight $\e\_{kt}$ of some arc $(k,t)$ is increased or a new $(k,t)$
arc is added to $\G$, and $\D p\_{ij},\;i,j\in V(\G),$ are the resulting increments of the accessibilities.
Then{\rm:}\\
{\rm (1)} $\D p\_{kt}>0;$\\
{\rm (2)} If $t$ mediates between $k$\/ and $i$, then $\D p\_{ki}>\D p\_{ti};$
if $k$ mediates between $i$\/ and $t$ then $\D p\_{it}>\D p\_{ik};$\\
{\rm (3)} (A)~If $t$ mediates between $k$\/ and $i$, then $\D p\_{kt}>\D
p\_{ki};$\\
\phantom{{\rm (3)}} (B)~If $k$ mediates between $i$\/ and $t$, then $\D p\_{kt}>\D p\_{it}$. }

\conditiont{Convexity} {(A)~If $p\_{ki}>p\_{ti}$ and $i\ne k,$ then there exists a $k$ to $i$ path such
that the difference $p\_{kj}-p\_{tj}$ strictly decreases as $j$ advances from $k$ to $i$ along this path.
(B)~If $p\_{ik}>p\_{it}$ and $i\ne k,$ then there exists an $i$ to $k$ path such that the difference
$p\_{jk}-p\_{jt}$ strictly increases as $j$ advances from $i$ to $k$ along this path.}

The results of testing $P^{{\rm out}}_{\tau}$ and $P^{{\rm in}}_{\tau}$ are collected{\x} in

\begin{theorem}
\label{otledostup}
The measures $P^{{\rm out}}_{\tau}$ and $P^{{\rm in}}_{\tau}$ satisfy
all the above conditions not partitioned into\/ $(A)$
                                          and\/~$(B)$.
Furthermore$,$
$P^{{\rm out}}_{\tau}$ obeys all\/~$(A)$ conditions and
$P^{{\rm in}}_{\tau}$        all\/~$(B)$ conditions.
\end{theorem}

\begin{proof}
{Let $P=P^{{\rm out}}_{\tau}.$
{\em Nonnegativity\/} follows from the definition of $P^{{\rm out}}_{\tau}$
and the positivity of arc weights.
{\em Reachability condition\/} follows from Proposition~\ref{reach}.
Item~(A) of the {\em self-accessibility condition\/} is true because, by
(\ref{Qtau}),
$p\_{ii}=\a\e(\FF^{i\rto i}(\tau)),$
$p\_{ij}=\a\e(\FF^{i\rto j}(\tau)),$
and
$\FF^{i\rto j}(\tau)\subset\FF^{i\rto i}(\tau)$ (strictly), where
$\a=\e^{-1}(\FF^{\rto}(\tau))$ and
$\FF^{i\rto j}(\tau)$ is the
$\FF^{i\rto j}$ for the digraph $\G(\tau)$ obtained from $\G$ by the
multiplication of all arc weights by~$\tau$; the same with
$\FF^{\rto}(\tau).$ Item~(A) of the {\em transit property\/} is proved
similarly.

To prove (A) of {\em convexity}, rewrite item~2 of
Theorem~\ref{th4} in the form $I=J(\tau)\,(I+\tau L).$
Consider entries $ki$ and $ti$ of $J(\tau)\,(I+\tau L)$. Since $i\ne k$
by assumption, we get
$p\_{ki}=\tau\sum_{j\ne i}\e\_{ji}(p\_{kj}-p\_{ki})$ and
$p\_{ti}=\tau\sum_{j\ne i}\e\_{ji}(p\_{tj}-p\_{ti})+\delta\_{it},$
where $\delta\_{it}=1$ if $i=t$ and $\delta\_{it}=0$ otherwise. Hence,
$p\_{ki}-p\_{ti}+\delta\_{it}
=\tau\sum_{j\ne i}\e\_{ji}\bigl((p\_{kj}-p\_{tj})-(p\_{ki}-p\_{ti})\bigr).$
Since $p\_{ki}-p\_{ti}>0$, there exists $j\in V(\G)$ such that
$\e\_{ji}\ne 0$ (and thus $(j,i)\in E(\G)$) and
$p\_{kj}-p\_{tj}>p\_{ki}-p\_{ti}$. Applying the same argument to
$j$ instead of $i$, and so forth, we finally obtain $k$ as the terminal
vertex of this path, as desired.

The {\em triangle inequality\/} follows from (A)~of {\em convexity\/} (taking $j=k$).
For the proof of {\em monotonicity\/} (items 1, 2, and~3A) we refer
to~\cite[Proposition~11]{AgaChe01}.

The corresponding statements for $P^{{\rm in}}_{\tau}$ follow similarly
or by duality.}
\end{proof}

It will be shown elsewhere that for a sufficiently small positive~$\tau,$
$P^{{\rm out}}_{\tau}$ additionally satisfies (B)~conditions, whereas
$P^{{\rm in}}_{\tau}$               satisfies (A)~conditions, and they
both satisfy the following {\it addition to monotonicity}:
{Suppose that the weight $\e\_{kt}$ of some arc $(k,t)$
increases or a new $(k,t)$ arc is added; then
for any $i,j\in V(\G),$ ($i\ne j$ or $k\ne t$) implies $\D
p\_{kt}>\D p\_{ij}.$
}

Consider now the accessibility measures
$\widetilde P^{{\rm out}}=(p\_{ij})=\vj
                        =\lim_{\tau\to\infty}P^{{\rm out}}_{\tau}\/$ and
$\widetilde P^{{\rm in}}=\lim_{\tau\to\infty}P^{{\rm  in}}_{\tau}$. Having in mind Theorem~\ref{M}, we call
$\q_{ij}$ the {\it limiting out-accessibility of $j$ from $i$}.

Let us say that a condition {\em is satisfied in the nonstrict form\/} if it is not generally satisfied, but
it becomes true after the substitution of $\ge$ for $>,$ $\le$ for $<$ and ``nonstrictly'' for ``strictly''
in the conclusion of this condition.

Similarly to Proposition~\ref{pr_dual} we have

\begin{proposition}
\label{pr_dual1}
The accessibility measures $\widetilde P^{{\rm out}}$ and
                           $\widetilde P^{{\rm in}}$ are dual.
\end{proposition}

The results of testing $\widetilde P^{{\rm out}}$ and $\widetilde P^{{\rm in}}$ are collected{\x} in

\begin{theorem}
\label{bliz}
The accessibility measures $\widetilde P^{{\rm out}}$ and
                           $\widetilde P^{{\rm in}}$
satisfy nonnegativity and the\/ {\rm ``$\Leftarrow$''} part of
reachability condition$,$ but they violate the\/ {\rm ``$\Rightarrow$''} part of reachability condition.
Moreover$,$
$\widetilde P^{{\rm out}}$ satisfies$,$ in the nonstrict
form$,$ items~$(A)$ of self-accessibility condition$,$ transit
property$,$ monotonicity$,$ and convexity$,$ whereas
$\widetilde P^{{\rm in}}$ satisfies in the nonstrict form items~$(B)$ of these conditions.
$\widetilde P^{{\rm out}}$ satisfies $(A)$ and
$\widetilde P^{{\rm in}}$  satisfies $(B)$ of triangle inequality for
proximities.
\end{theorem}

By virtue of Theorem~\ref{bliz}, the limiting accessibility measures only ``marginally'' correspond to the
conception of accessibility that underlies the above conditions.

\begin{proof}
{The nonstrict satisfaction of the conditions listed in the theorem follows from
Theorem~\ref{otledostup}, Proposition~\ref{pr_dual1} and item~5 of Theorem~\ref{th4}. To prove that the
strict forms of these conditions and the ``$\Rightarrow$'' part of reachability condition are violated, it
suffices to consider the digraph $\G$ with $n\ge3$, $E(\G)=\{(1,2),(2,3)\}$, and $\e\_{12}=\e\_{23}=1$.
}
\end{proof}

Let us mention one more class of accessibility measures, $(I+\a\J)^{-1}$, $0<\a<\si\_{(\di)}/\si\_{(\di+1)}$.
These measures are ``intermediate'' between $P^{{\rm out}}_{\tau}$ and $\widetilde P^{{\rm out}}$, because
they are positive linear combination of $J\_{(\di)}$ and $J\_{(\di+1)}$~\cite{AgaChe01}. That is why we
termed them the {\em matrices of dense out-forests}. In the terminology of \cite[p.~152]{MeyerStadelmaier78},
$(I+\a\J)^{-1}$ with various sufficiently small $\a>0$ make up a class of {\em nonnegative nonsingular
commuting weak inverses\/} for~$L$. These measures and the dual measures have been studied in \cite{AgaChe01}
(see also \cite[p.~270--271]{CheAga02a+}). Other interesting related topics are the forest distances
\cite{chebotarev02forest} and the forest based centrality measures~\cite{CheSha97}.

\section{Rooted forests and the problem of leaders}
\label{leader}

Ranking from tournaments or irregular pairwise contests is an old, but still intriguing problem. Its
statistical version is ranking objects on the base of paired comparisons~\cite{David88}. Analogous problems
of the analysis of individual and collective preferences arise in the contexts of policy, economics,
management science, sociology, psychology, etc.
Hundreds of methods have been proposed for handling these problems (for a review, see, e.g.,
\cite{David88,DavidAndrews93,CookKress92,BelkinLevin90,CheSha97a,CheSha99,Laslier97}).

In this section, we consider a weighted digraph~$\G$ that represents a competition (which need not be a round
robin tournament, i.e., can be ``incomplete'') with weighted pairwise results. The digraph can also represent
an arbitrary weighted preference relation. The result we present below can be easily extended to
multidigraphs.

One of the popular exquisite methods for assigning scores to the participants in a tournament was
independently proposed by Daniels \cite{Daniels69}, Moon and Pullman~\cite{MoonPullman69,MoonPullman70}, and
Ushakov~\cite{Ushakov71,Ushakov76} and reduces to finding nonzero and nonnegative solutions to the system of
equations
\begin{equation}
Lx=0.
\label{D}
\end{equation}

Entry $x_i$ of a solution vector $x=(x\_1\cdc x\_n)$ is considered as a sophisticated ``score'' attached to
vertex~$i$. This method was multiply rediscovered with different motivations (some references are given
in~\cite{CheSha99}). As Berman \cite{Berman80} noticed (although, in other contexts, similar results had been
obtained by Maxwell \cite{Maxwell1892} and other writers, see~\cite{CaplanZeilberger82}), if a digraph is
strong, then the general solution to (\ref{D}) is provided by the vectors proportional to $t=(t\_1\cdc
t\_n)^{\interca},$ where $t\_j$ is the weight of the set of spanning trees (out-arborescences) diverging
from~$j$. This fact can be easily proved as follows. By the matrix-tree theorem for digraphs (see,
e.g.,~\cite{Harary69}), $t\_j$ is the cofactor of any entry in the $j$th column of~$L$. Then for every $i\in
V(\G),$ $\sum^n_{j=1}\,\l\_{ij}\,t_j=\det L$ (the row expansion of $\det L$) and, since $\det L=0,\;$ $t$ is
a solution to~(\ref{D}). As $\rank(L)=n-1$ (since the cofactors of $L$ are nonzero), any solution to
(\ref{D}) is proportional to~$t$.

Berman \cite{Berman80} and Berman and Liu~\cite{BermanLiu96} asserted that this result is sufficient to rank
the players in an arbitrary competition, since the strong components of the corresponding digraph supposedly
``can be ranked such that every player in a component of higher rank defeats every player in a component of
lower rank. Now by ranking the players in each component we obtain a ranking of all the players.'' While the
statement about the existence of a natural order of the strong components is correct in the case of
round-robin tournaments, it need not be true for arbitrary digraphs that may have, for instance, several
source knots. That is why, the solution devised for strong digraphs does not enable one to rank the vertices
of an arbitrary digraph.

Let us consider the problem of interpreting, in terms of forests, the general solution to (\ref{D}) and the
problem of choosing a particular solution that could serve as a reasonable score vector in the case of
arbitrary digraph~$\G$.

If $\G$ contains more than one source knot, there is no spanning diverging tree in~$\G$. Recall that
$K\_1\cdc K\_{\di}$ are the source knots of $\G,$ where $\di$ is the out-forest dimension of $\G,$ and
$\ktil=\bigcup_{s=1}^{d'}K\_s$.

Suppose, without loss of generality, that the vertices of $\G$ are numbered as follows. The smallest numbers
are attached to the vertices in $K\_1$, the following numbers to the vertices in $K\_2$, etc., and the
largest numbers to the vertices in $V(\G)\setminus\ktil.$ Such a numeration we call {\em standard}.

\begin{theorem}
\label{d_sol}
Any column of $\J$ is a solution to~$(\ref{D})$.
Suppose that the numeration of vertices is standard and
$j\_1\in K\_1\cdc j\_{\di}\in K\_{\di}$. Then the columns
$\J\_{\bullet j\_1}\cdc\J\_{\bullet j\_{\di}}$ of\/ $\G$ make up an
orthogonal basis in the space of solutions to~$(\ref{D})$ and
$\J\_{\bullet j\_s}
=\e^{-1}(\TT\_s)
\bigl(0\cdc 0,\e(\TT_s^{i\_s+1})\cdc$
        $\e(\TT_s^{i\_s+k\_s}),0\cdc 0\bigr)^{\interca},$
where $\{i\_s+1\cdc i\_s+k\_s\}=K\_s$ and $\TT_s$ is the set of
out-arborescences of $K\_s,$ $s=1\cdc \di.$
\end{theorem}

By virtue of Theorem~\ref{d_sol}, the general solution to (\ref{D}) is the set of all linear combinations of
partial solutions that correspond to each source knot of~$\G$.

\begin{proof}
{The first statement follows from $L\q=0$ (item~3 of Theorem~\ref{th4}). By item~6 of Theorem~\ref{th4},
$\rank(\q)=\di$ and $\rank(L)=n-\di$. Hence, $\di$ is the dimension of the space of solutions to~(\ref{D}).
Let $j\_s\in K\_s,\;s=1\cdc\di.$ Then, by items~1 and~2 of Theorem~\ref{th2},
$$
\J\_{\bullet j\_s}
=\e^{-1}(\TT\_s)
\bigl(0\cdc 0,\e(\TT_s^{i\_s+1})\cdc\e(\TT_s^{i\_s+k\_s}),
      0\cdc 0\bigr)^{\interca}.
$$
These $\di$ solutions to~(\ref{D}) are orthogonal and thus, linearly independent.
}
\end{proof}

As a reasonable ultimate score vector, the arithmetic mean $x={1\over n}\J\cdot(1\cdc 1)^{\interca}$ of the
columns of $\q$ can be considered. A nice interpretation of this vector is given by

\begin{corollary}
\label{c-lim_distr}
{\rm (from Theorem~\ref{M}).}
For any Markov chain inversely corresponding to $\G,$ $x={1\over n}\J\cdot(1\cdc 1)^{\interca}$ is the limiting state
distribution$,$ provided that the initial state distribution is uniform.
\end{corollary}

It can be mentioned, however, that the ranking  method based on $\J$ takes into account long paths in $\G$
only. That is why, in any solutions to (\ref{D}), the vertices that are not in the source knots are assigned
zero scores, which is questionable. The estimates based on the matrices $Q(\tau),$ instead of $\q,$ are free
of this feature. On the other hand, both methods violate the {\em self-consistent monotonicity\/}
axiom~\cite{CheSha99}, and so do the methods that count the {\it routes\/} between vertices. This axiom is
satisfied by the {\it generalized Borda method\/} \cite{Che89,Che94} that produces the score vectors
$J'(\tau)\!\cdot\!(\od(1)-\id(1)\cdc\od(n)-\id(n))^{\interca}$, where $J'(\tau)$ is the matrix $J(\tau)$ of
the undirected graph corresponding to $\G$~\cite{Sha94}. In our opinion, the latter method can be recommended
as a well-grounded approach to scoring objects on the base of arbitrary weighted preference relations,
incomplete tournaments, irregular pairwise contests, etc.

\section*{A concluding remark: a communicatory interpretation of some forest matrices}

In closing, let us mention an interpretation of forest matrices in terms of information dissemination.
Consider the following metaphorical model. First, a plan of information transmission along a digraph is
chosen. Such a plan is a diverging forest $F\in\FF^{\rto}$: the information is injected into the roots of
$F$; then it ought to come to the other vertices along the arcs of $F$. Suppose that $\e\_{ij}\in ]0,1]$ is
the probability of successful information transmission along the $(i,j)$ arc, $i,j\in V(\G),$ and that the
transmission processes in different arcs are statistically independent. Then $\e(F)$ is the probability that
plan $F$ is successfully realized. Suppose now that each plan is selected with the same probability
$\card{\FF^{\rto}}^{-1}.$ Then $J\_{ij}$ (see~(\ref{J})) is the probability that the information came to $j$
from root $i$, provided that the transmission was successful. As a result, if one knows that the information
was corrupted at root $i$ and the transmission was successful, then $J\_{ij}$ is the probability that this
corrupted information came to~$j.$

Similarly, interpretations of this kind can be given to other normalized forest matrices. This model is
compatible with that of centered partitions \cite{Lenart98} and comparable with some models
of~\cite{Pavlov00}.

\bibliographystyle{endm}
\bibliography{c:/pavel/bibli/pavel/all2}

\begin{thebibliography}{10}
\expandafter\ifx\csname url\endcsname\relax
  \def\url#1{\texttt{#1}}\fi
\expandafter\ifx\csname urlprefix\endcsname\relax\def\urlprefix{URL }\fi
\newcommand{\enquote}[1]{``#1''}

\bibitem{AgaChe05LAA}
Agaev, R. and P.~Chebotarev, \emph{On the spectra of nonsymmetric {Laplacian}
  matrices}, Linear Algebra and Its Applications \textbf{399} (2005),
  pp.~157--168.

\bibitem{AgaChe00}
Agaev, R.~P. and P.~Y. Chebotarev, \emph{The matrix of maximum out forests 
  of a digraph and its applications}, Automation and Remote Control \textbf{61} 
  (2000), pp.~1424--1450.

\bibitem{AgaChe01}
Agaev, R.~P. and P.~Y. Chebotarev, \emph{Spanning forests of a digraph and
  their applications}, Automation and Remote Control \textbf{62} (2001),
  pp.~443--466.

\bibitem{BelkinLevin90}
Belkin, A.~R. and M.~S. Levin, \enquote{Prinyatie reshenii$:$ kombinatornye
  modeli approksimatsii informatsii $($Decision Making$:$ Combinatorial Models
  of Information Approximation$)$,} Nauka, Moscow, 1990, in Russian.

\bibitem{Berman80}
Berman, K.~A., \emph{A graph theoretical approach to handicap ranking of
  tournaments and paired comparisons}, SIAM Journal Alg. Disc. Meth. \textbf{1}
  (1980), pp.~359--361.

\bibitem{BermanLiu96}
Berman, K.~A. and Y.~Liu, \emph{Generalized bicycles}, Discrete Applied
  Mathematics \textbf{78} (1997), pp.~27--40.

\bibitem{CaplanZeilberger82}
Caplan, S.~R. and D.~Zeilberger, \emph{{T.L.~Hill}'s graphical method for
  solving linear equations}, Advances in Applied Mathematics \textbf{3} (1982),
  pp.~377--383.

\bibitem{chebotarev02spanning}
Chebotarev, P., \emph{Spanning forests of digraphs and limiting probabilities
  of {Markov} chains}, Electronic Notes in Discrete Mathematics \textbf{11}
  (2002), pp.~108--116.

\bibitem{Che04MCTT}
Chebotarev, P., \emph{Forest representations of {Markov} chain probabilities}
  (2004), unpublished.

\bibitem{CheAga02a+}
Chebotarev, P. and R.~Agaev, \emph{Forest matrices around the {Laplacian}
  matrix}, Linear Algebra and Its Applications \textbf{356} (2002),
  pp.~253--274.

\bibitem{chebotarev02forest}
Chebotarev, P. and E.~Shamis, \emph{The forest metrics for graph vertices},
  Electronic Notes in Discrete Mathematics \textbf{11} (2002), pp.~98--107.

\bibitem{Che89}
Chebotarev, P.~Y., \emph{Generalization of the row sum method for incomplete
  paired comparisons}, Automation and Remote Control \textbf{50} (1989),
  pp.~1103--1113.

\bibitem{Che94}
Chebotarev, P.~Y., \emph{Aggregation of preferences by the generalized row sum
  method}, Mathematical Social Sciences \textbf{27} (1994), pp.~293--320.

\bibitem{CheSha97a}
Chebotarev, P.~Y. and E.~V. Shamis, \emph{Constructing an objective function
  for aggregating incomplete preferences}, in: A.~Tangian and J.~Gruber,
  editors, \emph{Econometric Decision Models, Lecture Notes in Economics and
  Mathematical Systems}, Springer, Berlin, 1997 pp. 100--124.

\bibitem{CheSha97}
Chebotarev, P.~Y. and E.~V. Shamis, \emph{The matrix-forest theorem and
  measuring relations in small social groups}, Automation and Remote Control
  \textbf{58} (1997), pp.~1505--1514.

\bibitem{CheSha98a}
Chebotarev, P.~Y. and E.~V. Shamis, \emph{On a duality between metrics and
  ${\rm\Sigma}$-proximities}, Automation and Remote Control \textbf{59} (1998),
  pp.~608--612.

\bibitem{CheSha98}
Chebotarev, P.~Y. and E.~V. Shamis, \emph{On proximity measures for graph
  vertices}, Automation and Remote Control \textbf{59} (1998), pp.~1443--1459.

\bibitem{CheSha99}
Chebotarev, P.~Y. and E.~V. Shamis, \emph{Preference fusion when the number of
  alternatives exceeds two: Indirect scoring procedures}, Journal Franklin
  Inst. \textbf{336} (1999), pp.~205--226.

\bibitem{CookKress92}
Cook, W.~D. and M.~Kress, \enquote{Ordinal Information and Preference
  Structures$:$ Decision Models and Applications,} Prentice-Hall, Englewood
  Cliffs, 1992.

\bibitem{Daniels69}
Daniels, H.~E., \emph{Round-robin tournament scores}, Biometrika \textbf{56}
  (1969), pp.~295--299.

\bibitem{David88}
David, H.~A., \enquote{The Method of Paired Comparisons,} Griffin, London,
  1988, 2 edition.

\bibitem{DavidAndrews93}
David, H.~A. and D.~M. Andrews, \emph{Nonparametric methods of ranking from
  paired comparisons}, in: M.~A. Fligner and J.~S. Verducci, editors,
  \emph{Probability Models and Statistical Analyses for Ranking Data},
  Springer, New York, 1993 pp. 20--36.

\bibitem{FiedlerSedlacek58}
Fiedler, M. and J.~Sedl\'{a}\v{c}ek, \emph{O ${W}$-bas\'{\i}ch
  orientovan\'{y}ch graf\r{u}}, \v{C}asopis P\v{e}st. Mat. \textbf{83} (1958),
  pp.~214--225.

\bibitem{Harary69}
Harary, F., \enquote{Graph Theory,} Addison-Wesley, Reading, Mass., 1969.

\bibitem{Laslier97}
Laslier, J.~F., \enquote{Tournament Solutions and Majority Voting,} Springer,
  Berlin, 1997.

\bibitem{LeightonRivest83}
Leighton, T. and R.~L. Rivest, \emph{The {Markov} chain tree theorem}, Computer
  Science Technical Report MIT/LCS/TM-249, Laboratory of Computer Science, MIT,
  Cambridge, Mass. (1983).

\bibitem{LeightonRivest86}
Leighton, T. and R.~L. Rivest, \emph{Estimating a probability using finite
  memory}, IEEE Transactions on Information Theory \textbf{32} (1986),
  pp.~733--742.

\bibitem{Lenart98}
Lenart, C., \emph{A generalized distance in graphs and centered partitions},
  SIAM J. Discrete Math. \textbf{11} (1998), pp.~293--304.

\bibitem{Maxwell1892}
Maxwell, J.~C., \enquote{Electricity and Magnetism,} Oxford Univ. Press,
  London/New York, 1892, 3 edition, vol.~1, part.~{II}, p.~409.

\bibitem{MeyerStadelmaier78}
Meyer{,}~Jr., C.~D. and M.~W. Stadelmaier, \emph{Singular {M}-matrices and
  inverse positivity}, Linear Algebra and its Applications \textbf{22} (1978),
  pp.~139--156.

\bibitem{MoonPullman69}
Moon, J.~W. and N.~J. Pullman, \emph{Tournaments and handicaps}, in:
  \emph{Information Processing 68 $($Proc. IFIP Congress, Edinburgh, 1968$)$}
  (1969), pp. 219--223, vol.~1: Mathematics, Software.

\bibitem{MoonPullman70}
Moon, J.~W. and N.~J. Pullman, \emph{On generalized tournament matrices}, SIAM
  Review \textbf{12} (1970), pp.~384--399.

\bibitem{Pavlov00}
Pavlov, Y.~L., \enquote{Random Forests,} VSP, Utrecht, 2000.

\bibitem{Rothblum76a}
Rothblum, U.~G., \emph{Computation of the eigenprojection of a nonnegative
  matrix at its spectral radius}, Mathematical Programming Study \textbf{6}
  (1976), pp.~188--201.

\bibitem{Roubens96FSS}
Roubens, M., \emph{Choice procedures in fuzzy multicriteria decision analysis
  based on pairwise comparisons}, Fuzzy Sets and Systems \textbf{84} (1996),
  pp.~135--142.

\bibitem{Schwartz86}
Schwartz, T., \enquote{The Logic of Collective Choice,} Columbia Univ. Press,
  New York, 1986.

\bibitem{Sha94}
Shamis, E.~V., \emph{Graph-theoretic interpretation of the generalized row sum
  method}, Mathematical Social Sciences \textbf{27} (1994), pp.~321--333.

\bibitem{Ushakov71}
Ushakov, I.~A., \emph{The problem of choosing the preferable object}, Izv.
  Akad. Nauk SSSR. Tekhn. Kibernet. \textbf{No.~4} (1971), pp.~3--7, in
  Russian.

\bibitem{Ushakov76}
Ushakov, I.~A., \emph{The problem of choosing the preferred element: An
  application to sport games}, in: R.~E. Machol, S.~P. Ladany and D.~G.
  Morrison, editors, \emph{Management Science in Sports}, North-Holland,
  Amsterdam, 1976 pp. 153--161.

\bibitem{Vol'skii88}
Vol'skii, V.~I., \emph{Choice of best alternatives on directed graphs and
  tournaments}, Automation and Remote Control \textbf{49} (1988), pp.~267--278.

\bibitem{Zykov69}
Zykov, A.~A., \enquote{Teoriya Konechnykh Grafov $($Theory of Finite
  Graphs$)$,} Nauka, Novosibirsk, 1969, in Russian.

\end{thebibliography}

\end{document}